\documentclass[a4paper,12pt]{article}
\usepackage{amsmath}
\usepackage{latexsym}
\usepackage{url}
\usepackage{color}
\numberwithin{equation}{section}

\def\R{{\bf R}}

\def\N{{\bf N}}

\def\d{\displaystyle}
\def\e{{\varepsilon}}

\newtheorem{thm}{Theorem}[section]

\newtheorem{lem}{Lemma}[section]

\title{The lifespan of solutions of semilinear
wave equations with the scale-invariant damping\\
 in one space
dimension}
\author{
Masakazu Kato
\footnote{College of Liberal Arts, Mathematical Science Research Unit, 
Muroran Institute of Technology, 27-1, Mizumoto-cho, 
Muroran, Hokkaido 050-8585, Japan.
email: mkato@mmm.muroran-it.ac.jp.},
Hiroyuki Takamura
\footnote{Mathematical Institute,
Tohoku University,
 Aoba, Sendai 980-8578, Japan.
e-mail: hiroyuki.takamura.a1@tohoku.ac.jp.},
Kyouhei Wakasa
\footnote{Department of Creative Engineering, National Institute of Technology, Kushiro College, 2-32-1 Otanoshike-Nishi, Kushiro-Shi, Hokkaido 084-0916, Japan. e-mail: wakasa@kushiro-ct.ac.jp.
}
}
\date{
\[
\begin{array}{ll}
\mbox{\footnotesize{\bf Keywords:}}
& \mbox{\footnotesize semilinear wave equation, scale-invariant damping, lifespan}\\
\mbox{\footnotesize{\bf MSC2010:}}
& \mbox{\footnotesize primary 35L71, secondary 35B44}\\
\end{array}
\]
}
\begin{document}
\maketitle
\begin{abstract}
The critical constant $\mu$ (see (\ref{IVP1})) of time-decaying damping in the scale-invariant case 
is recently conjectured.
It also has been expected that the lifespan estimate is
the same as for the associated semilinear heat equations if the constant is in the \lq\lq heat-like" domain.
In this paper, we 
point out that this is not true if the total integral of the sum of initial position and speed vanishes.
In such a case, we have a new type of the lifespan estimates
which is closely related to the non-damped case in shifted space dimensions.
\end{abstract}


\section{Introduction}
\par
We consider the following initial value problem
for semilinear wave equations with the scale-invariant damping:
\begin{equation}\label{IVP1}
\left\{
\begin{array}{ll}
	\d v_{tt}-\Delta v+\frac{\mu}{1+t}v_t=|v|^p
	&\mbox{in}\quad \R^n\times[0,\infty),\\
	v(x,0)=\e f(x),\ v_t(x,0)=\e g(x),
	& x\in\R^n,
\end{array}
\right.
\end{equation}
where $p>1$, $\mu>0$,
the initial data $(f,g)\in H^1(\R^n)\times L^2(\R^n)$ is of compact support
and $\e>0$ is \lq\lq small".
The classification of general damping terms for the linear equation
is introduced by Wirth \cite{Wir1, Wir2, Wir3}.
The scale-invariant case is critical in the behavior of the solution.
For the outline of semilinear equations in other cases,
see Introduction of  Lai and Takamura \cite{LT}.
\par
It is interesting to look for the critical exponent $p_c(n)$ such that
\[
\left\{
\begin{array}{lll}
p>p_c(n)\ \mbox{(and may have an upper bound)}& \Longrightarrow & T(\e)=\infty,\\
1<p\le p_c(n) & \Longrightarrow & T(\e)<\infty,
\end{array}
\right.
\]
where $T(\e)$ is, the so-called lifespan, the maximal existence time of the energy solution of (\ref{IVP1})
with arbitrary fixed non-zero data.
Then we have the following conjecture:
\begin{equation}
\label{conj}
\left\{
\begin{array}{llll}
\mu\ge\mu_0(n) & \Longrightarrow & p_c(n)=p_F(n) & \mbox{(heat-like)},\\
0<\mu<\mu_0(n) & \Longrightarrow & p_c(n)=p_S(n+\mu) & \mbox{(wave-like)},
\end{array}
\right.
\end{equation}
where
\begin{equation*}
\mu_0(n):=\frac{n^2+n+2}{n+2}.
\end{equation*}
Moreover
\begin{equation*}
p_F(n):=1+\frac{2}{n}
\end{equation*}
is the so-called Fujita exponent which is the critical exponent of the  associated
semilinear heat equations $v_t-\Delta v=v^p$ and
\begin{equation*}
p_S(n):=\frac{n+1+\sqrt{n^2+10n-7}}{2(n-1)}\ (n\neq1),\ :=\infty\ (n=1)
\end{equation*}
is the so-called Strauss exponent which is the critical exponent of the associated
semilinear wave equations $v_{tt}-\Delta v=|v|^p$.
We note that $p_S(n)\ (n\neq1)$ is a positive root of
\begin{equation*}
\gamma(p,n):=2+(n+1)p-(n-1)p^2=0.
\end{equation*}
Moreover, $0<\mu<\mu_0(n)$ is equivalent to $p_F(n)<p_S(n+\mu)$.
Concerning the conjecture (\ref{conj}), D'Abbicco \cite{DABI} has obtained
heat-like existence partially as
\[
\mu\ge
\left\{
\begin{array}{cl}
5/3 & \mbox{for }n=1,\\
3 & \mbox{for }n=2,\\
n+2 & \mbox{for }n\ge3,
\end{array}
\right.
\]
while Wakasugi \cite{WY14_scale} has obtained blow-up
for $1<p\le p_F(n)$ and $\mu\ge1$, or $1<p\le p_F(n+\mu-1)$ and $0<\mu<1$.
We note that his result is the first blow-up result for super-Fujita exponents.
\par
Making use of the so-called Liouville transform
\[
u(x,t)=(1+t)^{\mu/2}v(x,t),
\]
one can rewrite (\ref{IVP1}) as
\begin{equation}
\label{IVP2}
\left\{
\begin{array}{ll}
\d u_{tt}-\Delta u+\frac{\mu(2-\mu)}{4(1+t)^2}u=\frac{|u|^p}{(1+t)^{\mu(p-1)/2}}
&\mbox{in}\quad \R^n\times[0,\infty),\\
u(x,0)=\e f(x),\ u_t(x,0)=\e \{\mu f(x)/2+g(x)\},
& x\in\R^n.
\end{array}
\right.
\end{equation}
Due to this observation,
D'Abbicco, Lucente and Reissig \cite{DLR15} have proved the wave-like part of
the conjecture (\ref{conj}) for $n=2,3$ when $\mu=2$.
We note that the radial symmetry is assumed for $n=3$ in \cite{DLR15}.
Moreover D'Abbicco and Lucente \cite{DL} have obtained the wave-like existence part
of (\ref{conj}) for odd $n\ge5$ when $\mu=2$ also with radial symmetry.
In the case $\mu=2$, (\ref{IVP2}) is a Cauchy problem for semilinear wave equations
with time-dependent coefficient on the right-hand side.
So, the regularity of the solution can be chosen higher,
sometimes a classical solution is handled.
For $\mu\neq2$, Lai, Takamura and Wakasa \cite{LTW} have first studied
the wave-like blow-up of the conjecture (\ref{conj}) with a loss
replacing $\mu$ by $\mu/2$ in the sub-critical case.
Initiating this, Ikeda and Sobajima \cite{IS} have obtained
the blow-up part of (\ref{conj}).
\par
For the semilinear wave equations with scale-invariant damping and mass,
the global existence of small data and the blow-up behavior were studied in
\cite{NPR17} and \cite{PR19}.
\par
For the lifespan estimate, one may expect that
\begin{equation}
\label{heat-lifespan}
T(\e)\sim
\left\{
\begin{array}{ll}
C\e^{-(p-1)/\{2-n(p-1)\}} & \mbox{for }1<p<p_F(n)\\
\exp\left(C\e^{-(p-1)}\right) & \mbox{for }p=p_F(n)
\end{array}
\right.
\end{equation}
for the heat-like domain $\mu\ge\mu_0(n)$ and
\begin{equation}
\label{wave-lifespan}
T(\e)\sim
\left\{
\begin{array}{ll}
C\e^{-2p(p-1)/\gamma(p,n+\mu)} & \mbox{for }1<p<p_S(n+\mu)\\
\exp\left(C\e^{-p(p-1)}\right) & \mbox{for }p=p_S(n+\mu)
\end{array}
\right.
\end{equation}
for the wave-like domain $0<\mu<\mu_0(n)$.
Here $T(\e)\sim A(\e,C)$ stands for the fact that there are positive constants,
$C_1$ and $C_2$, independent of $\e$ satisfying $A(\e,C_1)\le T(\e)\le A(\e,C_2)$.
Actually, (\ref{heat-lifespan}) for $n=1$ and $\mu=2$ is obtained by Wakasa \cite{wak16},
and (\ref{wave-lifespan}) is obtained by Kato and Sakuraba \cite{KS} for $n=3$ and $\mu=2$.
Also see Lai \cite{Lai} for the existence part of weaker solution.
Moreover, the upper bound of (\ref{heat-lifespan}) in the sub-critical case
is obtained by Wakasugi \cite{WY14_scale}.
Also the upper bound of (\ref{wave-lifespan}) is obtained by
Ikeda and Sobajima \cite{IS} in the critical case,
later it is reproved by Tu and Li \cite{TL2},
and Tu and Li \cite{TL1} in the sub-critical case.
\par
But we have the following fact. For the non-damped case, $\mu=0$, it is known that
(\ref{wave-lifespan}) is true for $n\ge3$, or $p>2$ and $n=2$.
The open part around this is $p=p_S(n)$ for $n\ge9$.
Other cases, (\ref{wave-lifespan}) is still true if the total integral of the initial speed vanishes,
i.e. $\int_{\R^n}g(x)dx=0$.
On the other hand,
we have 
\begin{equation}
\label{wave}
T(\e)\sim
\left\{
\begin{array}{ll}
C\e^{-(p-1)/2} & \mbox{for }n=1,\\
C\e^{-(p-1)/(3-p)} & \mbox{for $n=2$ and $1<p<2$},\\
Ca(\e) & \mbox{for $n=2$ and $p=2$}
\end{array}
\right.
\end{equation}
if $\int_{\R^n}g(x)dx\neq0$,
where $a=a(\e)$ is a positive number satisfying $\e^2a^2\log(1+a)=1$.
We note that (\ref{wave}) is smaller than 
the first line in (\ref{wave-lifespan}) with $\mu=0$ in each case.
For all the references of the case of $\mu=0$,
see Introduction of Imai, Kato, Takamura and Wakasa \cite{IKTW17}.
\par
Our aim in this paper is to show that
the lifespan estimates for (\ref{IVP2}) are similar to the ones for the non-damped case
even if $\mu$ is in the heat-like domain
by studying the special case of $n=1$ and $\mu=2\ge\mu_0(1)=4/3$.
That is, the result on (\ref{heat-lifespan}) by Wakasa \cite{wak16} mentioned above is true
only if $\int_{\R}\{f(x)+g(x)\}dx\neq0$.
More precisely, we shall show that
\begin{equation}
\label{1d_lifespan}
T(\e)\sim\left\{
\begin{array}{ll}
C\e^{-2p(p-1)/\gamma(p,3)} &\mbox{for} \quad\d 1<p<2,\\
Cb(\e) &\mbox{for} \quad\d p=2,\\
C\e^{-p(p-1)/(3-p)} &\mbox{for} \quad\d 2<p<3,\\
\exp(C\e^{-p(p-1)}) & \mbox{for}\quad\d p=p_F(1)=3
\end{array}
\right.
\end{equation}
if $\int_{\R}\{f(x)+g(x)\}dx=0$,
where $b = b(\e)$ is a positive number satisfying 
\begin{equation}
\label{b}
\e^2b\log(1+b)=1.
\end{equation} 
We note that (\ref{1d_lifespan}) is bigger than (\ref{heat-lifespan})
with $n=1$ and $\mu=2$ in each case.
This kind of phenomenon is observed also in two space dimensions
for $1<p\le p_F(2)=p_S(2+2)=2$ and $\mu=\mu_0(2)=2$.
Such a result will appear in our forthcoming paper \cite{IKTW}.
\par
This paper is organized as follows.
In the next section, we place precise statements on (\ref{1d_lifespan}).
Section 3, or 4, are devoted to the proof of the lower, or upper, bound of the lifespan respectively.


\section{Theorems and preliminaries}
\label{section:Pre}
\par
We shall show (\ref{1d_lifespan}) by establishing the following two theorems.

\begin{thm}
\label{T.1.1}
Let $n=1$, $\mu=2$ and $1<p\le 3=p_F(1)$.
Assume that $(f,g)\in C_0^2(\R)\times C_0^1(\R)$ satisfies
$\int_{\R}\{f(x)+g(x)\}dx=0$ and
\begin{equation}
\label{supp}
{\rm supp}\ (f,g) \subset\{x\in\R\ :\ |x|\le k\},\ k>1.
\end{equation}
Then, there exists a positive constant $\e_0=\e_0(f,g,p,k)$ such that a classical solution 
$u\in C^2(\R\times[0,T))$ of (\ref{IVP2}) exists as far as
\begin{equation}
\label{lower_lifespan}
T\le\left\{
\begin{array}{ll}
c\e^{-2p(p-1)/\gamma(p,3)} &\mbox{if} \quad\d 1<p<2,\\
cb(\e) &\mbox{if} \quad\d p=2,\\
c\e^{-p(p-1)/(3-p)} &\mbox{if} \quad\d 2<p<3,\\
\exp(c\e^{-p(p-1)}) & \mbox{if}\quad\d p=3
\end{array}
\right.
\end{equation}
for $0<\e\le\e_0$, where $c$ is a positive constant independent of $\e$
and $b(\e)$ is defined in (\ref{b}).
\end{thm}

\begin{thm}
\label{T.1.2}
Let $n=1$, $\mu=2$ and $1<p\le 3=p_F(1)$.
Assume that $(f,g)\in C_0^2(\R)\times C_0^1(\R)$ satisfy
$f(x)\ge0$ $(\not\equiv0)$, $f(x)+g(x)\equiv0$ and (\ref{supp}).
Then, there exists a positive constant $\e_1=\e_1(f,g,p,k)$ such that a classical solution 
$u\in C^2(\R\times[0,T))$ of (\ref{IVP2}) cannot exist whenever $T$ satisfies
\begin{equation*}
T\ge\left\{
\begin{array}{ll}
C\e^{-2p(p-1)/\gamma(p,3)} &\mbox{if} \quad\d 1<p<2,\\
Cb(\e) &\mbox{if} \quad\d p=2,\\
C\e^{-p(p-1)/(3-p)} &\mbox{if} \quad\d 2<p<3,\\
\exp(C\e^{-p(p-1)}) & \mbox{if}\quad\d p=3
\end{array}
\right.
\end{equation*}
for $0<\e\le\e_1$, where $C$ is a positive constant independent of $\e$
and $b(\e)$ is defined in (\ref{b}).
\end{thm}

As preliminaries for the proofs of the above theorems
we list some known facts. 
First, $u^0$ is defined by
\begin{equation}
\label{u^0}
u^0(x,t):=\frac{1}{2}\{f(x+t)+f(x-t)\}+\frac{1}{2}\int_{x-t}^{x+t}\{f(y)+g(y)\}dy
\end{equation}
with $(f,g)\in C^2(\R)\times C^1(\R)$ satisfies 
\[
\left\{
\begin{array}{ll}
\d u^0_{tt}-u^0_{xx}=0
& \mbox{in}\ \R\times[0,\infty),\\
u^0(x,0)=f(x),\ u^0_t(x,0)=f(x)+g(x),
& x\in \R.
\end{array}
\right.
\]
If we assume (\ref{supp}) and
\[
\int_{\R}\{f(x)+g(x)\}dx=0,
\]
then we have
\begin{equation}
\label{supp_u_L}
{\rm supp}\ u^0\subset\{(x,t)\in\R\times[0,\infty)\ :\ t-k\le |x|\le t+k\}.
\end{equation}
Moreover, if $u\in C(\R\times[0,\infty))$ is a solution of 
\begin{equation}
\label{IE}
u(x,t)=\e u^0(x,t)+L(|u|^p)(x,t)\quad\mbox{for}\ (x,t)\in\R\times[0,\infty),
\end{equation}
where
\begin{equation}
\label{L}
L(F)(x,t):=\frac{1}{2}\int_{0}^{t}\int_{x-t+s}^{x+t-s}\frac{F(y,s)}{(1+s)^{p-1}}dyds
\end{equation}
for $F\in C(\R\times[0,\infty))$,
then $u\in C^2(\R\times[0,\infty))$ is the solution to the initial value problem (\ref{IVP2}).
We also note that (\ref{supp}) implies
\begin{equation}
\label{supp_u}
{\rm supp}\ u\subset\{(x,t)\in\R\times[0,\infty)\ :\ |x|\le t+k\}.
\end{equation}
\par
We define the $L^{\infty}$ norm of $V$ by 
\begin{equation}\label{norm0}
\|V\|_0:=\sup_{(x,t)\in\R\times[0,T]}|V(x,t)|.
\end{equation}
Let $r=|x|$. 
For $r,t\ge0$, we define the following weighted functions:
\begin{equation}\label{weight}
w(r,t):=
	\left\{
	\begin{array}{lll}
		1 & \mbox{if}\ p>2,\\
		\{\log \tau_{+}(r,t)\}^{-1} & \mbox{if}\ p=2,\\
		\tau_{+}(r,t)^{p-2} & \mbox{if}\ 1<p<2,
\end{array}
\right.
\end{equation}
where we set
\begin{align*}
	\tau_+(r,t):=\frac{t+r+2k}{k}.
\end{align*}
For these weighted functions, we denote a weighted $L^\infty$ norm of $V$ by 
\begin{equation}\label{norm1}
	\|V\|:=\sup_{(x,t)\in\R\times[0,T]}\{w(|x|,t)|V(x,t)|\}.
\end{equation}
\par
Finally, we shall show some useful representations for $L$.
It is trivial that $1+s\ge (2k+s)/2k$ is valid for $s\ge0$ and $k>1$.
Setting $s=(\alpha+\beta)/2\ge0$ with $\alpha\ge0, \beta\ge-k$,
we have
\[
1+s\ge\frac{\alpha+2k}{4k},
\ \mbox{or}\ \ge\frac{\beta+2k}{4k}. 
\]
Thus, for $0 \leq \theta \leq 1$, we get 
\begin{equation}\label{est:weight}
	\frac{1}{1+s}\le \frac{4}{\{(\alpha+2k)/k\}^{\theta} 
	 \{(\beta+2k)/k\}^{1-\theta}}.
\end{equation}
Let $F=F(|x|,t)\in C(\R \times [0,T])$ and
\[
\mbox{supp}\ F\subset\{(x,t)\in \R\times[0,T] \ : \ |x|\le t+k\}. 
\]
From (\ref{L}), we obtain
\begin{equation*}
\begin{array}{ll}
|L(F)(x,t)|
&\d\le\frac{1}{2}\int_{0}^{t}ds\int_{r-t+s}^{r+t-s}\frac{|F(|y|,s)|}{(1+s)^{p-1}}dy\\
&=:L_{1}(F)(r,t)+L_{2}(F)(r,t),
\end{array}
\end{equation*}
where
\[
L_1(F)(r,t):=\frac{1}{2}\int_{0}^{t}ds\int_{|r-t+s|}^{r+t-s}\frac{|F(|y|,s)|}{(1+s)^{p-1}}dy
\]
and
\[
\begin{array}{ll}
L_2(F)(r,t)&:=
\d\frac{1}{2}\int_{0}^{(t-r)_+}\!\!\!\!ds\int_{r-t+s}^{t-r-s}\frac{|F(|y|,s)|}{(1+s)^{p-1}}dy\\
&\d=\int_{0}^{(t-r)_+}\!\!\!\!ds\int_{0}^{t-r-s}\frac{|F(|y|,s)|}{(1+s)^{p-1}}dy.
\end{array}
\]
Here we write $(a)_+=\max(a,0)$ for $a\in\R$. Changing the variables by $\alpha=s+y$, $\beta=s-y$ 
and making use of (\ref{est:weight}), we have
\begin{equation}
\label{dai2.3}
\begin{array}{l}
L_{1}(F)(r,t)\\
\d\le\int_{-k}^{t-r}d\beta
\int_{|t-r|}^{t+r}\frac{4^{p-2}|F((\alpha-\beta)/2,(\alpha+\beta)/2)|}
{ \{(\alpha+2k)/k\}^{\theta(p-1)}  \{(\beta+2k)/k\}^{(1-\theta)(p-1)}}d\alpha\\
\d\le\int_{-k}^{t+r}d\beta
\int_{\beta}^{t+r}\frac{4^{p-2}|F((\alpha-\beta)/2,(\alpha+\beta)/2)|}
{\{(\alpha+2k)/k\}^{\theta(p-1)}  \{(\beta+2k)/k\}^{(1-\theta)(p-1)}}d\alpha.
\end{array}
\end{equation}
Similarly it follows from (\ref{est:weight}) that
\begin{equation}
\label{dai2.4}
\begin{array}{l}
L_{2}(F)(r,t)\\
\d\le\int_{-k}^{t-r}d\beta
\int_{|\beta|}^{t-r}\frac{2^{-1}4^{p-1}|F((\alpha-\beta)/2,(\alpha+\beta)/2)|}	
{ \{(\alpha+2k)/k\}^{\theta(p-1)}  \{(\beta+2k)/k\}^{(1-\theta)(p-1)}}d\alpha\\
 \d\le\int_{-k}^{t+r}d\beta \int_{\beta}^{t+r}
 \frac{2^{-1}4^{p-1}|F((\alpha-\beta)/2,(\alpha+\beta)/2)|}
 { \{(\alpha+2k)/k\}^{\theta(p-1)} \{(\beta+2k)/k\}^{(1-\theta)(p-1)}}d\alpha.
\end{array}
\end{equation}
Therefore, we obtain by (\ref{dai2.3}) and (\ref{dai2.4}) that
\begin{equation}
\label{dai2.5}
\begin{array}{l}
|L(F)(x,t)|\\
\d\le\int_{-k}^{t+r}d\beta \int_{\beta}^{t+r}
 \frac{4^{p-1}|F((\alpha-\beta)/2,(\alpha+\beta)/2)|}
 { \{(\alpha+2k)/k\}^{\theta(p-1)} \{(\beta+2k)/k\}^{(1-\theta)(p-1)}}d\alpha.
\end{array}
\end{equation}

\section{Proof of Theorem\ref{T.1.1}}

First of all, we prove an estimate for the linear part of the solution from (\ref{IE}).
\begin{lem}\label{lm:linearest}
Let $u^0$ be as in {\rm (\ref{u^0})}.
Assume that the assumptions in Theorem \ref{T.1.1} are fulfilled. 
Then, there exists a positive constant $C_0$ such that 
\begin{equation}\label{linearest}
	\|u^0\|_0\le C_0.
\end{equation}
\end{lem}
\par\noindent
{\bf Proof.} It follows from (\ref{u^0}) and (\ref{supp_u_L}) that
\[
|u^0(x,t)|\le \|f\|_{L^{\infty}(\R)}+\|f+g\|_{L^{1}(\R)}.
\] 
Therefore, due to (\ref{norm0}), we obtain (\ref{linearest}). This completes the proof.
\hfill$\Box$
\vskip10pt

Next, we prove an a-priori estimate for a linear integral operator related
with the right-hand side of (\ref{IVP2}).
\begin{lem}\label{lm:apriori_u0}
Let $L$ be the linear integral operator defined by {\rm (\ref{L})}.
Assume that $V_0 \in C(\R\times[0,T])$
with {\rm  supp} $V_0\subset\{(x,t)\in \R\times[0,T] : t-k\le |x|\le t+k\}$. 
Then, there exists a positive constant $C_1$ independent of $T$ and $k$ such that 
\begin{equation}\label{apriori_u0}
\|L(|V_0|^p)\| \le C_1k^2\|V_0\|_0^{p}.
\end{equation}
\end{lem}
\par\noindent
{\bf Proof.}
We note that (\ref{apriori_u0}) follows from the following basic estimates:
\begin{equation}\label{basic.est1}
|L(\chi_{t-k\le r\le t+k})(x,t)|\le C_1k^2w(r,t)^{-1},
\end{equation}
where $\chi_A$ is a characteristic function of a set $A$.
\par
From now on to the end of this section,
$C$ stands for a positive constant independent of $T$ and $k$,
and may change from line to line.
It is easy to show (\ref{basic.est1}) by (\ref{dai2.5}) with $\theta=1$ and (\ref{weight}).
Actually we have that
\begin{align}
	|L(\chi_{t-k\le r\le t+k})(x,t)| &\le C\int_{-k}^{k}d\beta
	\int_{-k}^{t+r}\frac{d\alpha}{\{(\alpha+2k)/k\}^{p-1}}\nonumber\\
	&\leq Ck^2  \times
	\left\{
	\begin{array}{lll}
		1 & \mbox{if} &\d p>2,\\
		\d \log \tau_{+}(r,t)& \mbox{if} &\d p=2,\\
		\d \tau_{+}(r,t)^{2-p} & \mbox{if} &\d 1<p<2
	\end{array}
	\right.\nonumber\\
	&\le Ck^2 w(r,t)^{-1}.\nonumber
\end{align}
This completes the proof.
\hfill$\Box$
\vskip10pt

The following lemma contains one of the most essential estimates.
\begin{lem}\label{lm:apriori1}
Let $L$ be the linear integral operator defined by {\rm (\ref{L})}.
Assume that $V\in C(\R\times[0,T])$ with
{\rm  supp} $V\subset\{(x,t)\in \R\times[0,T] : |x|\le t+k\}$.
Then, there exists a positive constant $C_2$ independent of $T$ such that 
\begin{equation}\label{apriori1}
	\|L(|V|^p)\|\le C_2k^2\|V\|^pD(T),
\end{equation}
where $D(T)$ is defined by 
\begin{equation}\label{D}
D(T):=
\left\{
\begin{array}{lll}
\log T_k & \mbox{if}\quad p=3,\\
T_k^{3-p} & \mbox{if}\quad 2<p<3,\\
 T_k \log T_k & \mbox{if}\quad p=2,\\
T_k^{\gamma(p,3)/2} & \mbox{if}\quad 1<p<2
\end{array}
\right.
\end{equation}
with $T_k:=(T+2k)/k$.
\end{lem}
\par\noindent
{\bf Proof.}
We note that (\ref{apriori1}) follows from the following basic estimates:
\begin{equation*}
	|L(w^{-p})(x,t)|\le C_2k^2D(T)w(r,t)^{-1}. 
\end{equation*}
We divide the proof into three cases.
\par
\noindent (i) \ Case of $2<p \leq 3$.
\par
It follows from (\ref{weight}), (\ref{dai2.5}) with $\theta=1$ and (\ref{D}) that
\begin{align}
	|L(w^{-p})(x,t)| &\le \d C\int_{-k}^{t+r}d\beta
	\int_{\beta}^{t+r}\frac{d\alpha}{\{(\alpha+2k)/k\}^{p-1}}\nonumber\\
	&\le \d Ck\int_{-k}^{t+r}\{(\beta+2k)/k\}^{2-p}d\beta\nonumber\\
	&\leq Ck^2 \times
	\left\{
	\begin{array}{ll}
		\log \tau_{+}(r,t) & (p=3)\\
		\tau_{+}(r,t)^{3-p} & (2<p<3)
	\end{array}
	\right.\nonumber\\
	&\le Ck^2D(T)w(r,t)^{-1}.\nonumber
\end{align}
Here we have used by (\ref{supp_u}) that
\[
\tau_+(r,t)\le\frac{2t+3k}{k}\le 2T_k
\quad\mbox{and}\quad T_k\ge2.
\]
From now on, we will employ this estimate at the end of each case.
\par
\noindent (ii) \ Case of $p=2$.
\par
It follows from (\ref{dai2.5}) with $\theta=1/2$, (\ref{weight}) and (\ref{D}) that
\[
\begin{array}{l}
|L(w^{-p})(x,t)|\\
\d\le C\int_{-k}^{t+r}d\beta\int_{\beta}^{t+r} \frac{\log^2\{(\alpha+2k)/k\}}
{\left\{ (\alpha+2k)/k \right\}^{1/2} 
\left\{ (\beta+2k)/k \right\}^{1/2}} d\alpha\\
\d\le C\log^2\tau_+(r,t)
\int_{-k}^{t+r} \left(\frac{\beta+2k}{k}\right)^{-1/2} d\beta
\int_{-k}^{t+r} \left(\frac{\alpha+2k}{k}\right)^{-1/2}  d\alpha\\
\le Ck^2 \tau_{+}(r,t)\log^2\tau_{+}(r,t)\\
\le Ck^2D(T) w(r,t)^{-1}.
\end{array}
\]
\par
\noindent (iii) \ Case of $1<p<2$.
\par
Similarly to the above, it follows from (\ref{dai2.5}) with $\theta=1$ and  (\ref{weight})
that
\[
\begin{array}{llll}
	|L(w^{-p})(x,t)|&\le \d C\int_{-k}^{t+r}d\beta
	\int_{-k}^{t+r} \left(\frac{\alpha+2k}{k}\right)^{(2-p)p-(p-1)} d\alpha\nonumber\\
	&\leq Ck^2 \tau_{+}(r,t)^{-p^2+p+3}\nonumber\\
	&\le Ck^2D(T)w(r,t)^{-1}.
\end{array}
\]
The proof is now completed.
\hfill$\Box$
\vskip10pt

Finally, we state an a-priori estimate of mixed type.

\begin{lem}\label{lm:apriori3}
Let $L$ be the linear integral operator defined by {\rm (\ref{L})},
and $V,D(T)$ be as in Lemma \ref{lm:apriori1}.
Assume that $V_0 \in C(\R\times[0,T])$ with
\[
\mbox{\rm  supp}\ V_0\subset\{(x,t)\in \R\times[0,T]\ :\ t-k\le |x|\le t+k\}. 
\]
Then, there exists a positive constant $C_3$ independent of $T$ and $k$ such that 
\begin{equation*}
	\|L(|V_0|^{p-1}|V|)\| \le C_3 k^2 \|V_0\|_0^{p-1}\| V \| D(T)^{1/p}.
\end{equation*}
\end{lem}
{\bf Proof.}
Similarly to the proof of Lemma \ref{lm:apriori_u0}, we shall show 
\begin{equation}\label{basic.est3}
	|L(\chi_{t-k\le r\le t+k}w^{-1})(x,t)|\le C_3k^2w(r,t)^{-1}D(T)^{1/p}.
\end{equation}
\par
\noindent (i) \ Case of $2< p \leq 3$.
\par
Since $w(r,t)=1$,
(\ref{basic.est3}) is established by the estimates for $2<p\le3$ in Lemma \ref{lm:apriori1} 
and $1 \leq D(T)^{1/p}$. 
\par
\noindent (ii) \ Case of $p=2$.
\par
It follows from (\ref{dai2.5}) with $\theta=1$ and (\ref{weight}) that
\begin{align*}
	|L(\chi_{t-k\le r\le t+k}w^{-1})(x,t)|
	&\le\d C\int_{-k}^{k}d\beta
	\int_{-k}^{t+r}\frac{\log\{(\alpha+2k)/k\}}{(\alpha+2k)/k}d\alpha\\
	&\le Ck^2 \log^2\tau_{+}(r,t)\\
	&\le Ck^2 \log T_k\cdot w(r,t)^{-1}.
\end{align*}
Since $\log T_{k} \le D(T)^{1/2}$, we obtain (\ref{basic.est3}).
\par
\noindent (iii) \ Case of $1<p<2$.
\par
It follows from (\ref{dai2.5}) with $\theta=1$ that
\begin{align*}
	|L(\chi_{t-k\le r \le t+k}w^{-1})(x,t)| &\le\d C\int_{-k}^{k}d\beta
	\int_{-k}^{t+r} \left(\frac{\alpha+2k}{k}\right)^{2-p-(p-1)}d\alpha\\
	&\le Ck^2T_k^{2-p}w(r,t)^{-1}.
\end{align*}
Since $2-p\le \gamma(p,3)/2p$, we obtain (\ref{basic.est3}).\\
The proof is now completed.
\hfill$\Box$
\vskip10pt
\par
\noindent
{\bf Proof of Theorem \ref{T.1.1}.}
We consider the following integral equation.
\begin{align}
U=L(|\e u^{0}+U|^{p}) \quad \mbox{in} \quad \R \times [0, T].
\label{dai4.1}
\end{align}
Suppose we have a solution $U=U(x,t)$ of (\ref{dai4.1}).
Then, by putting $u=U+\e u^{0}$, we obtain a solution of (\ref{IE}) and its lifespan
is the same as that of $U$.
Thus, our aim here is to  construct a solution of (\ref{dai4.1}) in the Banach space,
\begin{align}
	X:=\{ U(x,t) \in C(\R \times [0,T]) \ : \ 
	 \mbox{supp} \ U
	\subset \{(x,t):|x| \leq t+k \} \}
	\nonumber
\end{align}
which is equipped with the norm (\ref{norm1}).
\par
Define a sequence of functions $\{ U_{l} \} \subset X$ by
\begin{align}
	U_{1}=0,\quad U_{l}=L(|\e u^0+U_{l-1}|^p)\quad \mbox{for} \quad 
	l \geq 2\nonumber
\end{align}
and set
\begin{align}
	&M_{0}:=2^{p-1}  C_{1} k^2 C_{0}^p,\nonumber\\
	&C_{4}:=(2^{2(p+1)}p)^p
	\max\{C_{2}k^2M_{0}^{p-1}, (C_{3}k^2C_{0}^{p-1})^{p}\},\nonumber
\end{align}
where $C_{i}$ $(0 \leq i \leq 3)$ are positive constants given 
in Lemma {\ref{lm:linearest}}, Lemma {\ref{lm:apriori_u0}}, Lemma {\ref{lm:apriori1}}
and Lemma {\ref{lm:apriori3}}.
Then, analogously to the proof of Theorem 1 in \cite{IKTW17},
we see that  $\{  U_{l} \}$ is a Cauchy sequence in $X$ provided
that  the inequality
\begin{align}
	C_{4}\e^{p(p-1)}D(T) \leq 1\label{dai4.2}
\end{align}
holds. Since $X$ is complete, there exists a function $U$ such that $U_{l}$ converges
to $U $ in $X$. Therefore $U$ satisfies (\ref{dai4.1}).
\par
Note that (\ref{lower_lifespan}) follows from (\ref{dai4.2}).
We shall show this fact only in the case of $p=2$ since the other cases can be proved similarly.
By definition of $b$ in (\ref{b}), we know that $b(\e)$ is decreasing in $\e$
and $\displaystyle \lim_{\e \to 0+0} b(\e)= \infty$.
Let us fix $\e_0>0$ as 
\begin{align}
	1<C_{5}b(\e_0),\label{dai4.3}
\end{align}
where $\displaystyle C_{5}=\min \left\{  2^{-1}, (3C_{4})^{-1} \right\}$.
For $0< \e\leq \e_0$, we take $T$ to satisfy
\begin{align}
	1\leq T<C_{5} b(\e).\label{dai4.4}
\end{align}
Since $k>1$, it follows from (\ref{D}) and  (\ref{dai4.4}) that
\begin{align*}
	C_{4} \e^2 D(T) 
	& \leq C_{4} \e^2 (3T) \log (2T+1)\\
	& \leq 3C_{4} C_{5} \e^2 b(\e) \log (2C_{5}b(\e)+1)\\
	& \leq b(\e) \e^2 \log(b(\e)+1)= 1.
\end{align*}
Hence, if we assume (\ref{dai4.3}) and (\ref{dai4.4}), then (\ref{dai4.2}) holds.
Therefore (\ref{lower_lifespan}) in the case $p=2$ is obtained for $0 < \e \leq \e_{0}$.
This completes the proof of Theorem \ref{T.1.1}.
\hfill$\Box$

\section{Proof of Theorem\ref{T.1.2}}
\par
In order to obtain an upper bound of the lifespan,
we shall take a look on the ordinary differential inequality for
\[
F(t):=\int_{\R}u(x,t)dx
\]
and shall follow the arguments in Section 5 of Takamura \cite{Takamura15}.
The equation in (\ref{IVP2}) with $\mu=2$ and (\ref{supp_u}) imply that
\begin{equation}
\label{int|u|^p}
F''(t)=\frac{1}{(1+t)^{p-1}}\int_{\R}|u(x,t)|^pdx
\quad\mbox{for}\ t\ge0.
\end{equation}
Hence, H\"older's inequality and (\ref{supp_u}) yield that
\begin{equation}
\label{F''}
F''(t)\ge2^{-(p-1)}(t+k)^{-2(p-1)}|F(t)|^p
\quad\mbox{for}\ t\ge0.
\end{equation}
Due to the assumption on the initial data in Theorem \ref{T.1.2},
\[
f(x)\ge0 \ (\not\equiv0),\quad f(x)+g(x)\equiv0,
\]
we have
\begin{equation}
\label{initial}
F(0)>0,\quad F'(0)=0.
\end{equation}
Neglecting the nonlinear term in (\ref{IE}), from (\ref{u^0}) and (\ref{supp}), we also obtain
the following point-wise estimate.
\begin{equation}
\label{point}
u(x,t)\ge\frac{1}{2}f(x-t)\e
\quad\mbox{for}\ x+t\ge k\ \mbox{and}\ -k\le x-t\le k.
\end{equation}
\par
First, we shall handle the sub-critical case.
In such a case,
the following basic lemma is useful.

\begin{lem}[\cite{Takamura15}]
\label{lem:Kato}
Let $p>1, a>0, q>0$ satisfy
\begin{equation}
\label{M}
M:=\frac{p-1}{2}a-\frac{q}{2}+1>0.
\end{equation}
Assume that $F\in C^2([0,T))$ satisfies
\begin{eqnarray}
& F(t)\ge  At^a &  \mbox{for}\ t\ge T_0,\label{ineq:A}\\
& F''(t)\ge  B(t+k)^{-q}|F(t)|^p & \mbox{for}\ t\ge0,\label{ineq:B}\\
& F(0)>0,\ F'(0)=0,\label{ineq:C}&
\end{eqnarray}
where $A,B,k,T_0$ are positive constants.
Moreover, assume that
there is a $t_0>0$ such that  
\begin{equation}
\label{t_0}
F(t_0)\ge 2F(0).
\end{equation}
Then, there exists a positive constant $C_*=C_*(p,a,q,B)$ such that
\begin{equation}
\label{est:T_1}
T<2^{2/M}T_1
\end{equation}
holds provided
\begin{equation}
\label{condi2}
T_1:=\max\left\{T_0,t_0,k\right\}\ge C_* A^{-(p-1)/(2M)}.
\end{equation}
\end{lem}
This is exactly Lemma 2.2 in \cite{Takamura15}, so that we shall omit the proof here.
We already have (\ref{ineq:B}) and (\ref{ineq:C}), so that
the key estimate is (\ref{ineq:A}) which is expected better than a constant $F(0)$
trivially follows from (\ref{ineq:B}).
\par
From now on to the end of this section,
$C$ stands for a positive constant independent of $\e$,
and may change from line to line.
It follows from (\ref{int|u|^p}) and (\ref{point}) that
\[
F''(t)\ge\frac{1}{(1+t)^{p-1}}\int_{t-k}^{t+k}|u(x,t)|^pdx
\ge C\e^pt^{1-p}
\quad\mbox{fot}\ t\ge k.
\]
Since (\ref{ineq:B}) and (\ref{ineq:C}) imply $F(t)>0$ and $F'(t) \geq 0$ for $t\geq0$, integrating this inequality twice in $t$, we obtain
\begin{equation}
\label{F}
F(t)\ge C\e^p\times
\left\{
\begin{array}{ll}
t^{3-p} & \mbox{if}\ 1<p<2,\\
\d t\log\frac{t}{2k} & \mbox{if}\ p=2,\\
t & \mbox{if}\ p>2
\end{array}
\right.
\quad\mbox{for}\ t\ge4k.
\end{equation}
\par\noindent
(i) Case of $1<p<2$.
\par
According to (\ref{F}), one can apply Lemma \ref{lem:Kato} to our situation with
\[
A=C\e^p,\ a=3-p>0,\ B=2^{-(p-1)},\ q=2(p-1).
\]
In this case, the blow-up condition (\ref{M}) is satisfied by
\[
2M=(p-1)(3-p)-2(p-1)+2=\frac{\gamma(p,3)}{2}>0.
\]
Next we fix $t_0$ to satisfy (\ref{t_0}).
Due to (\ref{F}), it is
\[
F(t_0)\ge C\e^pt_0^{3-p}=2F(0)= 2\|f\|_{L^1(\R)}\e,
\]
namely
\[
t_0=C\e^{-(p-1)/(3-p)}.
\]
Hence, setting
\[
T_0=C_*A^{-(p-1)/(2M)}=C\e^{-2p(p-1)/\gamma(p,3)},
\]
we have a fact that there exists an $\e_1=\e_1(f,g,p,k)>0$ such that
\[
T_1:=\max\{T_0,t_0,k\}=T_0=C\e^{-2p(p-1)/\gamma(p,3)}\ge4k
\]
holds for $0<\e\le\e_1$ because of
\[
\frac{1}{3-p}<\frac{2p}{\gamma(p,3)}
\quad\Longleftrightarrow\quad p>1.
\]
Therefore, from (\ref{est:T_1}), we obtain $T < 2^{2/M}T_1=C\e^{-2p(p-1)/\gamma(p,3)}$ as desired.\\
\par\noindent
(ii) Case of $2<p<3$.
\par
According to (\ref{F}), one can apply Lemma \ref{lem:Kato} to our situation with
\[
A=C\e^p,\ a=1,\ B=2^{-(p-1)},\ q=2(p-1).
\]
In this case, the blow-up condition (\ref{M}) is satisfied by
\[
2M=p-1-2(p-1)+2=3-p>0.
\]
Next we fix $t_0$ to satisfy (\ref{t_0}).
Due to (\ref{F}), it is
\[
F(t_0)\ge C\e^pt_0=2F(0)= 2\|f\|_{L^1(\R)}\e,
\]
namely
\[
t_0=C\e^{-(p-1)}.
\]
Hence, setting
\[
T_0=C_*A^{-(p-1)/(2M)}=C\e^{-p(p-1)/(3-p)},
\]
we have a fact that there exists an $\e_1=\e_1(f,g,p,k)>0$ such that
\[
T_1:=\max\{T_0,t_0,k\}=T_0=C\e^{-p(p-1)/(3-p)}\ge4k
\]
holds for $0<\e\le\e_1$ because of
\[
1<\frac{p}{3-p}
\quad\Longleftrightarrow\quad p>\frac{3}{2}.
\]
Therefore we obtain $T <  2^{2/M}T_1=C\e^{-p(p-1)/(3-p)}$ as desired.\\
\par\noindent
(iii) Case of $p=2$.
\par
Neglecting the logarithmic term in (\ref{F}),
similarly to the case of $2<p<3$,
one can apply Lemma \ref{lem:Kato} to our situation with
\[
A=C\e^2,\ a=1,\ B=2^{-1},\ q=2,\ 2M=1.
\]
We shall fix a $T_0$ as follows.
In order to establish (\ref{condi2}) in Lemma \ref{lem:Kato},
we have to assume that $T_0\ge C_*A^{-1}$ namely
\[
A\ge C_*T_0^{-1}.
\]
On the other hand, (\ref{ineq:A}) in Lemma \ref{lem:Kato}
can be established by (\ref{F}) as far as
\[
C\e^2\log\frac{T_0}{2k}\ge A.
\]Hence $T_0$ must satisfy
\begin{equation}
\label{T_0}
\e^2T_0\log\frac{T_0}{2k}\ge C_{**},
\end{equation}
where $C_{**}$ is a positive constant independent of $\e$.
Here we identify a constant $C$ as $C_{**}$ to fix $T_0$.
Recall the definition of $b(\e)$ in (\ref{b})
and the fact that $b(\e)$ is monotonously decreasing in $\e$ and
$\lim_{\e\rightarrow0+0}b(\e)=\infty$.
If $C_{**} \geq 1$, then we set $T_0=4kC_{**}b(\e)$. Taking $\e$ small to satisfy $C_{**}b(\e)\ge1$, we have
\[
\e^2T_0\log\frac{T_0}{2k}\ge4kC_{**}\e^2b(\e)\log\{1+C_{**}b(\e)\}\ge4kC_{**}.
\]
Therefore (\ref{T_0}) holds if $C_{**}\ge1$ by $k>1$.
On the other hand, if $C_{**}<1$, then we set $T_0=4kb(\e)$.
In this case, taking $\e$ small to satisfy $b(\e)\ge1$, we have
\[
\e^2T_0\log\frac{T_0}{2k}\ge4k\e^2b(\e)\log\{1+b(\e)\}=4k,
\]
so that (\ref{T_0}) holds by $4k>1>C_{**}$.
In this way one can say that our situation can be applicable to Lemma \ref{lem:Kato}
with $T_0=Cb(\e)$ for small $\e$ except for $t_0$ in (\ref{t_0}).
\par
In this case, (\ref{t_0}) follows from (\ref{F}) and
\[
F(t_0)\ge C\e^2t_0\log\frac{t_0}{2k}=2F(0)=2\|f\|_{L^1(\R)}\e,
\]
namely
\[
\e t_0\log\frac{t_0}{2k}=C.
\] 
Comparing this equality with (\ref{T_0}), we know that
there exists an $\e_1=\e_1(f,g,k)>0$ such that
\[
T_1:=\max\{T_0,t_0,k\}=T_0=Cb(\e)\ge4k
\]
holds for $0<\e\le\e_1$.
Therefore we obtain $T < 2^{2/M}T_1=Cb(\e)$ as desired.\\
\par\noindent
(iv) Case of $p=p_F(1)=3$
\par
Even in this case, (\ref{F}) is still valid.
But $a=1$ and $p=3$ yield $M=0$ in Lemma \ref{lem:Kato}.
So we need a critical version of the lemma,
which is a variant of Lemma 2.1 in Takamura and Wakasa \cite{TW11}
with a slightly different initial condition.
One can readily show it by small modification.
Here we shall avoid to employ it, and shall make use of (\ref{ineq:B}) and (\ref{F})
only to give a simple proof by means of \lq\lq slicing method" of the blow-up domain
introduced in Agemi, Kurokawa and Takamura \cite{AKT}.
\par
For $j\in\N\cup\{0\}$, define
\[
a_j:=\sum_{i=0}^j\frac{1}{2^i}\quad\mbox{and}\quad K:=4k.
\]
Assume presumably
\begin{equation}
\label{j}
F(t)\ge D_jt\log^{b_j}\frac{t}{a_jK}
\quad\mbox{for}\quad t\ge a_jK,
\end{equation}
where each $b_j$ and $D_j$ are positive constants.
We note that (\ref{j}) with $j=0$ is true by (\ref{F}) if we set $b_0=0$ and $D_0=C\e^3$.
Plugging (\ref{j}) into the right hand side of (\ref{F''}) with a restriction 
to the interval $[a_jK,\infty)$,
we obtain that
\[
F''(t)\ge 2^{-6}D_j^3t^{-1}\log^{3b_j}\frac{t}{a_jK}
\quad\mbox{for}\quad t\ge a_jK
\]
which yields that
\[
F'(t)\ge 2^{-6}D_j^3\cdot\frac{1}{3b_j+1}\log^{3b_j+1}\frac{t}{a_jK}
\quad\mbox{for}\quad t\ge a_jK.
\]
Integrating this inequality and diminishing the interval to make use of
\[
\int_{a_jK}^t\log^{3b_j+1}\frac{s}{a_jK}ds
\ge\int_{a_jt/a_{j+1}}^t\log^{3b_j+1}\frac{s}{a_jK}ds
\quad\mbox{for}\quad t\ge a_{j+1}K,
\]
we obtain that
\[
F(t)\ge 2^{-6}D_j^3\cdot\frac{1}{3b_j+1}\left(1-\frac{a_j}{a_{j+1}}\right)t\log^{3b_j+1}\frac{t}{a_{j+1}K}
\quad\mbox{for}\quad t\ge a_{j+1}K.
\]
Thus, due to
\[
1-\frac{a_j}{a_{j+1}}=\frac{1}{2^{j+1}a_{j+1}}\ge\frac{1}{2^{j+2}},
\]
(4.14) inductively holds if the sequence $\{b_j\}$ is defined by
\begin{equation}
\label{b_j}
b_{j+1}=3b_j+1,\ b_0=0\quad\mbox{for}\ j\in\N\cup\{0\}
\end{equation}
and $\{D_j\}$ is defined by\begin{equation}
\label{D_j}
D_{j+1}:=\frac{D_j^3}{2^{j+8}(3b_j+1)},\ D_0:=C\e^3\quad\mbox{for}\ j\in\N\cup\{0\}.
\end{equation}
It is easy to see that (\ref{b_j}) gives us  
\begin{equation}
\label{b_j:sol}
b_j=\frac{3^j-1}{2}\quad\mbox{for}\ j\in\N\cup\{0\}.
\end{equation}
\par
From now on, let us look for a suitable lower bound of $D_j$ by (\ref{D_j}).
Since
\[
3b_j+1=b_{j+1}\le\frac{3^{j+1}}{2}
\quad\mbox{for}\ j\in\N\cup\{0\}
\]
by (\ref{b_j:sol}), we have
\[
\log D_{j+1}\ge 3\log D_j-(2j+8)\log3\quad\mbox{for}\  j\in\N\cup\{0\}
\]
which yields
\[
\log D_j\ge3^{j-1}\log D_0-\log3\sum_{i=0}^{j-1}3^{j-1-i}(2i+8)\quad\mbox{for}\  j\in\N.
\]
Hence, it follows from
\[
S:=\lim_{j\rightarrow\infty}\sum_{i=0}^{j-1}\frac{2i+8}{3^i}>0
\]
by d'Alembert criterion that
\[
D_j\ge\left(\frac{D_0}{3^S}\right)^{3^{j-1}}\quad\mbox{for}\ j\in\N.
\]
Therefore, together with (\ref{j}), we have
\[
F(t)\ge\left(\frac{D_0}{3^S}\right)^{3^{j-1}}t\log^{(3^j-1)/2}\frac{t}{2K}
=\frac{3^S}{D_0}t \left( \log^{-1/2}\frac{t}{2K} \right) I(t)^{3^j}
\]
for $t\ge 2K$ and $j\ge1$, where we set
\[
I(t):=\frac{D_0}{3^S}\log^{1/2}\frac{t}{2K}.
\]
This inequality means that
\[
\lim_{j\rightarrow\infty}F(t_1)=\infty
\]
if there exists a $t_1\ge2K$ such that $I(t_1)>1$.
It can be achieved by
\[
\exp \left( -\left(\frac{D_0}{3^S}\right)^{-2} \right) \frac{t_1}{2K}>1.
\]
Therefore, $T$ has to satisfy 
\[
T\le2K\exp\left(C\e^{-6}\right).
\]
The proof is now completed in all the  cases.
\hfill$\Box$

\section*{Acknowledgement}
\par
This work started when the second author was working in Future University Hakodate
and third author was working in Muroran Institute of Technology.
The second author has been partially supported by
Special Research Expenses in FY2017, General Topics (No.B21), Future University Hakodate,
also by the Grant-in-Aid for Scientific Research (B) (No.18H01132) and (C) (No.15K04964), 
Japan Society for the Promotion of Science.


\bibliographystyle{plain}

\end{document}